\documentclass[12pt]{article}
\usepackage{amsmath}
\usepackage{amssymb}
\newsymbol \blackbox 1004
\newcommand{\eh}{\hfill}\newlength{\sperr}

\def\a{\alpha}
\def\b{\beta}
\def\g{\gamma}

\def\s{\sigma}
\def\S{\Sigma}

\def\t{\theta}

\def\wh{\widehat}
\def\wt{\widetilde}
\def\vt{\vartheta}
\def\BC{{\mathbb C}}
\def\BR{{\mathbb R}}

\def\im{{\rm Im\ }}

\newtheorem{Pa}{Paper}[section]
\newtheorem{Tm}[Pa]{{\bf Theorem}}
\newtheorem{La}[Pa]{{\bf Lemma}}
\newtheorem{Cy}[Pa]{{\bf Corollary}}
\newtheorem{Rk}[Pa]{{\bf Remark}}

\newtheorem{Ee}[Pa]{{\bf Example}}
\newtheorem{Dn}[Pa]{{\bf Definition}}
\newtheorem{Pn}[Pa]{{\bf Proposition}}
\newcommand{\CC}
{{\mathchoice {\setbox0=\hbox{$\displaystyle\rm C$}\hbox{\hbox
to0pt{\kern0.4\wd0\vrule height0.9\ht0\hss}\box0}}
{\setbox0=\hbox{$\textstyle\rm C$}\hbox{\hbox
to0pt{\kern0.4\wd0\vrule height0.9\ht0\hss}\box0}}
{\setbox0=\hbox{$\scriptstyle\rm C$}\hbox{\hbox
to0pt{\kern0.4\wd0\vrule height0.9\ht0\hss}\box0}}
{\setbox0=\hbox{$\scriptscriptstyle\rm C$}\hbox{\hbox
to0pt{\kern0.4\wd0\vrule height0.9\ht0\hss}\box0}}}}

\title{Discrete skew selfadjoint canonical systems and the isotropic
Heisenberg magnet model}

\author{\normalsize  M.A. KAASHOEK and A.L. SAKHNOVICH
}

\date{}
\begin{document}
\maketitle

\begin{abstract}
A discrete analog of a skew selfadjoint canonical (Zakharov-Shabat
or AKNS) system with a pseudo-exponential potential is introduced.
For the corresponding Weyl function the direct and inverse problem
are solved explicitly in terms of  three parameter matrices. As an
application explicit solutions are obtained for the discrete
integrable nonlinear equation corresponding to the isotropic
Heisenberg magnet model. State space techniques from mathematical
system theory play an important role in the proofs.
\end{abstract}

\setcounter{section}{-1}
\section{Introduction} \label{intro}
\setcounter{equation}{0}

 In this paper we shall treat a discrete
analog of the well-known skew selfadjoint canonical (Dirac type,
Zakharov-Shabat or AKNS) system:
\begin{equation} \label{0.1}
-i J  \frac{d Y}{dx}(x,z)=  z Y(x, z )+  V(x)Y(x, z ),\  J =
\left[
\begin{array}{cc}
I_{m} & 0 \\ 0 & -I_{m}
\end{array}
\right], \  x\geq 0.
\end{equation}
Here $z$ is a spectral variable, $Y$ and $V$ are  $2m \times 2m $
matrix functions on the half line, and $V$ is skew selfadjoint,
that is, $V(x)^*=JV(x)$ with  $V(x)^*$ being the matrix adjoint of
$V(x)$. To obtain the discrete analog of (\ref{0.1}) let $U$ be
the unique solution of the initial value problem
\begin{equation} \label{0.3}
 \frac{d U}{dx}(x)=-i U(x)JV(x), \quad x\geq 0,\quad U(0)=I_{2m}.
\end{equation}
Since $JV(x)$ is selfadjoint, we get from (\ref{0.3}) that $U(x)$
is unitary for each $x\geq 0$. Now put $S(x)=U(x)J U(x)^*$  and
$W(x, z) = U(x)Y(x,z)$. Then
\begin{equation} \label{0.4}
\frac{dW}{dx}(x, z)=i z S(x) W(x, z), \quad S(x)=S(x)^*=S(x)^{-1},
\quad x\geq 0.
\end{equation}
It is now immediate that
\begin{equation} \label{0.5}
W_{n+1}( \lambda)-W_n(\lambda)=- \frac{i}{\lambda}S_n
W_n(\lambda), \quad S_n=S_n^*=S_n^{-1},\quad n=0, \, 1, \ldots
\end{equation}
is a natural discrete analog of (\ref{0.1}). This discrete analog
of the continuous pseudo-canonical system is very important. In
fact, when $m=1$, then the system (\ref{0.5}) turns out to be an
auxiliary system for the nonlinear isotropic Heisenberg magnet
(IHM) model \cite{Skl} (see also the detailed discussion after
Theorem \ref{Tm2.3} below and the historical remarks in
\cite{FT}). Motivated by the IHM model we shall use the term
\emph{spin sequence} to denote any sequence of $N \times N$
matrices $\{S_n\}$ satisfying
\begin{equation} \label{0.5'}
S_n=S_n^*=S_n^{-1},\quad n=0, \, 1, \, 2, \ldots
\end{equation}

As for the skew  selfadjoint continuous case \cite{SaA1} (see also
\cite{GKS2}), one can associate with (\ref{0.5}) an $m \times m$
matrix function $\varphi (\lambda)$, meromorphic on $\Im \lambda
<- \delta <0$, such that
\begin{equation} \label{2.1}
\sum_{n=0}^{\infty}[ \varphi(\lambda)^* \quad I_m]\,W_n(\lambda)^*
W_n(\lambda) \left[
\begin{array}{c}
\varphi(\lambda)  \\ I_m
\end{array}
\right] < \infty,
\end{equation}
where $W_n(\lambda)$, $n\geq 0$, is the \emph{fundamental
solution} of (\ref{0.5}), i.e., the $2m\times 2m$ matrix solution
$W_n(\lambda)$ of (\ref{0.5}) normalized by the condition
$W_0(\lambda)=I_{2m}$. One refers to $\varphi$ as the \emph{Weyl
function} of $(\ref{0.5})$. When the Weyl function is rational and
strictly proper we shall recover the system (\ref{0.5}) explicitly
from its Weyl function. For this purpose we need to introduce spin
sequences that are the discrete analogs of the pseudo-exponential
potentials from \cite{GKS2,GKS6} (see also the references
therein).

The spin sequences from this  special class are defined in terms
of three parameter matrices with the following properties. First
fix an integer $N>0$, and consider an $N \times N$ matrix $\a$
with $\det \, \a \not= 0$, an  $N \times N$ matrix $\S_0$ such
that $\S_0=\S_0^*$, and  an $N \times 2m$ matrix $\Lambda_0$.
These matrices should satisfy the following matrix identity
\begin{equation} \label{0.8}
\a \S_0- \S_0 \a^*= i \Lambda_0 \Lambda_0^*.
\end{equation}
Given these three matrices $\a$, $\S_0$, and $\Lambda_0$ we define
for $n=1,2,\ldots$ the $N \times 2m$ matrix $\Lambda_n$ and the $N
\times N$ matrix $\S_n$ via recursion:
\begin{eqnarray}
 \Lambda_{n+1}&=& \Lambda_n+i \a^{-1} \Lambda_n J,\nonumber \\
 \noalign{\vskip6pt}
\S_{n+1}&=&\S_n+ \a^{-1} \S_n (\a^*)^{-1}+ \a^{-1} \Lambda_n J
\Lambda_n^* (\a^*)^{-1}.\label{0.7}
\end{eqnarray}
Next assume that  the matrices $\S_n$, $n=0,1,2,\ldots$, are
non-singular. Then we say that the sequence of matrices $\{ S_n
\}$ defined by
\begin{equation} \label{0.6}
S_n=J+ \Lambda_n^* \S_n^{-1} \Lambda_n - \Lambda_{n+1}^*
\S_{n+1}^{-1} \Lambda_{n+1}, \quad n=0,1,2,\ldots,
\end{equation}
is the \emph{spin sequence determined} by the parameter matrices
$\a$, $\S_0$ and $\Lambda_0$. Notice that this requires the
invertibility of the matrices $\S_n$.

For spin sequences defined in this way our first theorem presents
an formula for the fundamental solution $W_n(\lambda)$ of
(\ref{0.5}).

\begin{Tm} \label{Tm1.1}
Let  $\a$ $( \det \, \a \not= 0)$, $\S_0$ $(\S_0=\S_0^*)$ and
$\Lambda_0$ satisfy $(\ref{0.8})$, and
 assume that $\det \, \S_n \not=0$ for $0 \leq n \leq
M$, where $\S_n$ is given by $(\ref{0.7})$. For $0 \leq n \leq
M-1$ let $ S_n $ be the  matrices determined  by $\a$, $\S_0$ and
$\Lambda_0$ via $(\ref{0.6})$ and $(\ref{0.7})$. Then
$S_n=S_n^*=S_n^{-1}$ for $0 \leq n \leq M-1$, and for $0 \leq n
\leq M$ the fundamental solution $W_n(\lambda )$ of the discrete
system $(\ref{0.5})$ can be represented in the form
\begin{equation}\label{1.-1}
W_n( \lambda )=W_{ \alpha , \Lambda }(n, \lambda ) \Big(
I_{2m}-\frac{i}{\lambda}J \Big)^{n} W_{ \alpha , \Lambda }(0,
\lambda )^{-1},
\end{equation}
where $W_{ \alpha , \Lambda }(n, \lambda )$ is defined by
\begin{equation}\label{1.0}
W_{ \alpha , \Lambda }(n, \lambda )=I_{2m}+i \Lambda_n^{*}
\Sigma_n^{-1} ( \lambda I_{N}- \alpha )^{-1} \Lambda_n.
\end{equation}
\end{Tm}

When $\S_0>0$, there exist simple conditions on $\a$ and
$\Lambda_0$ to guarantee that $\det \, \S_n \not= 0$. First, if
$\S_0>0$, then without loss of generality we can assume that
$\S_0=I_N$. Indeed, it is easy to see that the sequence of
matrices $\{ S_n \}$ defined by (\ref{0.6}) and (\ref{0.7}) does
not change if we substitute $\a$, $\S_0$ and $\Lambda_0$ by
$\S_0^{-\frac{1}{2}} \a \S_0^{\frac{1}{2}}$, $I_N$  and
$\S_0^{-\frac{1}{2}} \Lambda_0$. So let us assume that $\S_0=I_N$.
Next, we partition  $\Lambda_0$ into two $N \times m$ blocks
$\t_1$ and $\t_2$ as follows: $\Lambda_0=[\t_1 \quad \t_2]$. This
together with $\S_0=I_N$ allows us to rewrite (\ref{0.8}) in the
form
\begin{equation} \label{2.0}
\a - \a^*=i(\t_1 \t_1^* + \t_2 \t_2^*).
\end{equation}
Furthermore, in this case $\Lambda_n$ is given by
\begin{equation} \label{2.14}
\Lambda_n=[(I_N+i \a^{-1})^n \t_1 \quad (I_N-i \a^{-1})^n \t_2].
\end{equation}
Finally, we shall assume that the pair  $\{\a$, $\t_1\}$ is
\emph{full range} which means that
\[
{\mathbb C}^N= {\mathrm{span}}\{\a^k \t_1\BC^m\mid k=0, 1,2, \dots
, N-1 \}.
\]
The following proposition shows that under these conditions
automatically $\det \, \a \not=0$ and $\det \, \S_n \not= 0$ for
$n=0,1,2,\dots$.

\begin{Pn} \label{Pn2.0} Let  $\a$ be a square matrix of order $N$,
and  $\t_1$ and $\t_2$ be $N\times m$ matrices satisfying
$(\ref{2.0})$. Assume that the pair $\{\a$, $\t_1\}$ is full
range. Then all the eigenvalues of $\a$ are in the open upper half
plane $\BC_+$, and for $n=1,2,\ldots$ the matrices $\S_n$ defined
by $(\ref{0.7})$, with $\S_0=I_N$ and $\Lambda_n$ given by
$(\ref{2.14})$, are positive definite and satisfy the identity
\begin{equation} \label{1.1}
\a \S_n- \S_n \a^*= i \Lambda_n \Lambda_n^*.
\end{equation}
\end{Pn}

\begin{Dn} \label{DnAdm}
A triple of matrices $\a$, $\t_1$ and $\t_2$, with $\a$ square of
order $N$ and $\t_1$ and $\t_2$ of size $N\times m$, is called
admissible if the pairs $\{\a$, $\t_1\}$ and $\{\a$, $\t_2\}$ are
full range and the identity $($\ref{2.0}$)$ holds.
\end{Dn}
We denote by the acronym FG (\emph{finitely generated}) the class
of spin sequences $\{S_n \}$ determined by the matrices $\a$,
$\S_0=I_N$ and $\Lambda_0=[\t_1 \quad \t_2]$, where $\a$, $\t_1$
and $\t_2$ form an admissible triple. In this case we also say
that these spin sequences are \emph{determined} by the
corresponding admissible triples. The next two theorems present
the solutions of the direct and inverse problem in terms of the
Weyl function.

\begin{Tm} \label{Tm2.1} Assume that the spin sequence
$\{ S_n \}_{n \geq 0}$ of the discrete pseudo-canonical system
$(\ref{0.5})$ belongs to the class FG and is determined by the
admissible triple $\a$, $\t_1$ and $\t_2$. Then the system
$(\ref{0.5})$ has a unique Weyl function $\varphi$, which
satisfies $(\ref{2.1})$ on the half plane $\Im \lambda <-
\frac{1}{2}$,
 a finite number of poles excluded,
 and this function is given by the formula
\begin{equation} \label{2.2}
\varphi ( \lambda )=i \theta_{1}^{*}( \lambda I_{N} - \beta )^{-1}
\theta_{2},
\end{equation}
where
\begin{equation} \label{2.3}
\beta = \alpha -i \theta_{2} \theta_{2}^{*}.
\end{equation}
 \end{Tm}

Notice that the function $\varphi$ in (\ref{2.2}) is a strictly
proper $m\times m$ rational matrix function. Conversely, if
$\varphi$  is a  strictly proper $m\times m$ rational matrix
function, then it admits a representation of the form
\begin{equation}  \label{2.24}
\varphi ( \lambda ) =i \vartheta^{*}_{1}( \lambda I_{n}- \gamma
)^{-1}
 \vartheta_{2},
\end{equation}
where $\gamma$ is a square matrix and
$\vartheta_{1},\vartheta_{2}$ are matrices of size $n\times m$. We
refer to the right hand side of (\ref{2.24}) as a \emph{minimal
realization} of $\varphi$ if among all possible representations
(\ref{2.24}) of $\varphi$ the order $n$ of the matrix $\gamma$ is
as small as possible. This terminology is taken from mathematical
system theory. We can now state the solution of the inverse
problem.

\begin{Tm} \label{Tm2.3}
Let $ \varphi $ be a strictly proper rational $m \times m$ matrix
function, given by the minimal realization $($\ref{2.24}$)$. There
is a unique  positive definite $n\times n$ matrix solution $X$ of
the algebraic Riccati equation
\begin{equation}  \label{2.23}
\gamma X - X
\gamma^{*}=i(X\vartheta_{1}\vartheta_{1}^*X-\vartheta_{2}\vartheta_{2}^*).
\end{equation}
Using $X$ define matrices
 $ \theta_{1}, \theta_{2}$, and $ \alpha =
 \beta +i \theta_{2} \theta_{2}^{*} $
 by
 \begin{equation}  \label{2.25}
\theta_{1}=X^{ \frac{1}{2}} \vartheta_{1}, \hspace{1em}
\theta_{2}=X^{- \frac{1}{2}} \vartheta_{2}, \hspace{1em} \beta
=X^{- \frac{1}{2}} \gamma X^{ \frac{1}{2}}.
\end{equation}
Then $ \alpha $, $ \theta_{1}$, and $\theta_{2}$ form an
admissible triple, and the given matrix function $\varphi$ is the
Weyl function of a system $(\ref{0.5})$ of which the spin sequence
$\{ S_n \} \in$ FG and is uniquely  determined by the admissible
triple $ \alpha $, $ \theta_{1}$, and $\theta_{2}$.
\end{Tm}

Next, we describe connections with the nonlinear IHM equation. For
this purpose consider the zero curvature representation \cite{Skl}
of the IHM model:
\begin{equation} \label{0.5.-2}
\frac{d }{d t}G_n(t, \lambda)=F_{n+1}(t, \lambda )G_n(t, \lambda
)-G_n(t, \lambda )F_n(t, \lambda),
\end{equation}
where
\begin{equation} \label{0.5.-1}
G_n(t, \lambda)=I_{2}- \frac{i}{\lambda}S_n(t), \quad F_n(t,
\lambda)= \frac{V_n^+(t)}{\lambda - i}+ \frac{V_n^-(t)}{\lambda +
i},
\end{equation}
\begin{equation} \label{0.5.4}
V_r^{\pm}(t):=(1+ \overrightarrow{S}_{r-1}(t) \, \cdot \,
\overrightarrow{S}_{r}(t))^{-1}(I_2 \pm S_r(t))(I_2 \pm
S_{r-1}(t)).
\end{equation}
Here the vectors $\overrightarrow{S}_{r}=[S^1_r \quad S^2_r \quad
S^3_r]$ belong to $\BR^3$, $\BR$ is real axis, the dot $\cdot$
denotes the scalar product in $\BR^3$, and the correspondence
between the spin matrix $S_r$ and the spin vector
$\overrightarrow{S}_{r}$ is given by the equality
\begin{equation} \label{0.5.2}
S_r= \left[
\begin{array}{lr}
S_r^3 & S_r^1 - i S_r^2 \\ S_r^1 + i S_r^2 & - S_r^3
\end{array}
\right].
\end{equation}
In other words the IHM equation
\begin{equation} \label{0.5.1}
\frac{d \overrightarrow{S}_n}{d t}=2 \overrightarrow{S}_n
\bigwedge \left( \frac{\overrightarrow{S}_{n+1}}{1+
\overrightarrow{S}_n \, \cdot \, \overrightarrow{S}_{n+1}} +
\frac{\overrightarrow{S}_{n-1}}{1+ \overrightarrow{S}_{n-1} \,
\cdot \, \overrightarrow{S}_{n}} \right),
\end{equation}
where $\bigwedge$ stands for the  vector product in $\BR^3$, is
equivalent to the compatibility condition (\ref{0.5.-2}) of the
systems
\begin{eqnarray}
W_{n+1}(t, \lambda)&=&G_n(t, \lambda)W_n(t,\lambda),\nonumber\\
\noalign{\vskip6pt} \frac{d}{dt}W_n(t, \lambda) &=&F_n(t,
\lambda)W_n(t, \lambda)\quad (n \geq 0).\label{0.5.0}
\end{eqnarray}
In (\ref{0.5.1}) it is required  that
$||\overrightarrow{S}_{n}||=1$. Now one can see easily  that the
representation (\ref{0.5.2}), where $\overrightarrow{S}_{n} \,
\cdot \overrightarrow{S}_{n}=1$, is equivalent to the equalities
(\ref{0.5'}) with $S_n \not= \pm I_2$. Thus  the first  system in
(\ref{0.5.0}) coincides with system
 (\ref{0.5}), where $m=2$, and $S_n
\not= \pm I_2$. We use these connections to obtain explicit
solutions of the IHM model.

The literature on continuous canonical systems is very rich,
especially for the selfadjoint case; see, for instance, the books
\cite{dB,FT,GK,LS,SaL3}. Selfadjoint continuous canonical systems
with pseudo-exponential potentials have been introduced in
\cite{GKS1}, and for this class of potentials various direct and
inverse problems have been solved; see \cite{GKS6} and the
references therein. The subclass of strictly pseudo-exponential
potentials has been treated in \cite{AG2,AGKS,AD}. Interesting
recent results on the spectral theory of selfadjoint discrete
systems and various useful references on this subject can be found
in \cite{AG1,CG,DKS,FG,RS,T}. Mainly  Jacobi matrices (or block
Jacobi matrices as in \cite{AG1}) that are related to  Toda chain
problems have been studied. For the skew selfadjoint discrete case
some references can be found in \cite{EJ,FT,H}.

Theorem \ref{Tm1.1} is the discrete analog of Theorem 1.2 in
\cite{GKS2}. The right hand side of (\ref{1.0}) can be viewed as
the transfer function of a linear input-output system (see
\cite{BGK}). Transfer functions of the special form given by
(\ref{1.0}) were introduced in \cite{SaL1}, and also used for the
representation of the fundamental solutions of continuous
canonical systems \cite{SaL2,SaL3}. In Theorem \ref{Tm1.1} we are
closer to \cite{GKS2} (see also \cite{SaALAA,SaA3}), where the
dependence on the parameter $x$ differs from the one in
\cite{SaL2,SaL3}.
 The condition on an admissible triple $\a$, $\t_1$,  $\t_2$
that the pairs $\{\a$, $\t_1 \}$ and $\{ \a$, $\t_2 \}$ are full
range pairs is specific for the discrete case. Nevertheless
Theorems \ref{Tm2.1}, \ref{Tm2.3} and parts of their proofs  are
analogous to results and proofs in Section 2 of \cite{GKS2}.

This paper consists of four sections not counting this
introduction. Since elements from mathematical system theory play
an important role in this paper, we present in the first section
the necessary preliminaries from that area. In the second section
we prove Theorem \ref{Tm1.1} and present some auxiliary results
that will be used in the application to the IHM model. Theorems
\ref{Tm2.1}, \ref{Tm2.3} and Proposition \ref{Pn2.0} are proved in
Section 3. In Section 4 we construct solutions of the IHM equation
(\ref{0.5.1}), describe the evolution of the Weyl function and
consider a simple example.
\section{Preliminaries from mathematical system \\ theory} \label{Ap}
\setcounter{equation}{0}

The material from the state space theory of rational matrix
functions, that is used in this paper, has its roots in the Kalman
theory of input-output systems \cite{KFA}, and can be found in
books, see, e.g.,  \cite{ka,CF}. In general the rational matrix
functions appearing in this paper are {\em proper}, that is,
analytic at infinity, and they are square, of size $m\times m$,
say. Such a function $F$ can be represented in the form
\begin{equation}
\label{app.1} F(\lambda)=D+C(\lambda I_N-A)^{-1}B,
\end{equation}
where $A$ is a square matrix (of which the order $N$  may be much
larger than m), the matrices $B$ and $C$ are of sizes $N\times m$
and $m\times N$, respectively, and $D=F(\infty)$. In this paper
$D$ is often a zero matrix, and in that case $F$ is called
\emph{strictly proper}. The representation (\ref{app.1}) is called
a {\em realization} or a {\em transfer matrix representation} of
$F$, and the number $N=\mathrm{ord}(A)$, is called the {\em state
space dimension} of the realization. Here $\mathrm{ord}(A)$
denotes the \emph{order} of the matrix $A$.

Realizations of a fixed $F$ are not unique. The realization
(\ref{app.1}) is said to be {\em minimal} if its state space
dimension $N$ is minimal among all possible realizations of $F$.
The state space dimension of a minimal realization of $F$ is
called the {\em McMillan degree} of $F$ and is denoted by $\deg
F$. Notice that $\deg F=0$ corresponds to the case when
$\mathrm{ord}(A)=0$, and this occurs if and only if
$F(\lambda)\equiv
D$. The realization (\ref{app.1}) of $F$ is minimal if and only if
\begin{equation}
\label{app.2'} {\mathrm{span}}\bigcup_{k=0}^{N-1}\im
A^kB=\BC^N,\quad \bigcap_{k=0}^{N-1}\mathrm{Ker}\, CA^k=\{0\},
\quad N= \mathrm{ord}(A).
\end{equation}
If for a pair of matrices $\{A$, $B \}$ the first equality in
(\ref{app.2'}) holds, then $\{A$, $B \}$ is called {\em
controllable} or a {\em full range pair}. If the second equality
in (\ref{app.2'}) is fulfilled, then $\{C$, $A \}$ is said to be
{\em observable} or a {\em zero kernel pair}. If a pair $\{A$, $B
\}$ is full range, and $K$ is an $m\times N$ matrix, where $N$ is
the order of $A$ and $m$ is the number of columns for $B$, then
the pair $\{A-BK$, $B \}$ is also full range. An analagous result
holds true for zero kernel pairs.

Minimal realizations are unique up to a basis transformation, that
is, if (\ref{app.1}) is a minimal realization of $F$ and
$F(\lambda)=D+\widetilde C(\lambda I_N -\widetilde
A)^{-1}\widetilde B$ is a second minimal realization of $F$, then
there exists an
invertible matrix $S$ such that
\[
 \widetilde A=SAS^{-1},\quad \widetilde B=SB,\quad
\widetilde C=CS^{-1}.
\]
In this case $S$ is called a {\em state space similarity}.

Finally if in (\ref{app.1}) we have $D=I_m$, then $F(\lambda)$ is
invertible whenever $\lambda$ is not an eigenvalue of $A-BC$ and
in that case
\begin{equation} \label{inv}
F(\lambda)^{-1}=I_m-C(\lambda I_N-A^\times)^{-1}B,\quad
A^\times=A-BC.
\end{equation}
\section{The fundamental solution} \label{fund}
\setcounter{equation}{0} In this section we prove Theorem
\ref{Tm1.1} and present some results that will be used in Section
\ref{IHM} (a result on the invertibility of matrices $\S_n$, in
particular).

\medskip\noindent
{\bf Proof of Theorem \ref{Tm1.1}.} First we shall show that
equalities (\ref{0.8}) and (\ref{0.7}) yield the identity
(\ref{1.1}) for all $n \geq 0$. The statement is proved by
induction. Indeed, for $n=0$ it is true by assumption. Suppose
(\ref{1.1}) is true for $n=r$. Then using the expression for
$\S_{r+1}$ from (\ref{0.7}) and identity (\ref{1.1}) for $n=r$ we
get
\begin{eqnarray}
\a \S_{r+1}- \S_{r+1} \a^*= \qquad  \qquad \qquad \qquad \qquad
\nonumber
\\ \quad = i \Lambda_r \Lambda_r^* +
 i \a^{-1} \Lambda_r \Lambda_r^* (\a^*)^{-1}+
 \Lambda_r J \Lambda_r^* (\a^*)^{-1}
 -
\a^{-1} \Lambda_r J \Lambda_r^*. \label{1.3}
\end{eqnarray}
The first relation in (\ref{0.7}) and formula (\ref{1.3}) yield
(\ref{1.1}) for $r=n+1$ and thus for all $n \geq 0$.

The next equality will be crucial for our proof. Namely we shall
show that for $0 \leq n \leq M-1$ we have
\begin{equation} \label{1.4}
W_{ \alpha , \Lambda }(n+1, \lambda )(I_{2m}-
\frac{i}{\lambda}J)=(I_{2m}- \frac{i}{\lambda}S_n)W_{ \alpha ,
\Lambda }(n, \lambda ).
\end{equation}
By (\ref{1.0}) formula (\ref{1.4}) is equivalent to the formula
\begin{eqnarray}
\frac{1}{\lambda}(S_n-J)&=&(I_{2m}- \frac{i}{\lambda}S_n)
\Lambda_n^{*} \Sigma_n^{-1} ( \lambda I_{N}- \alpha )^{-1}
\Lambda_n -\nonumber\\ \noalign{\vskip6pt} &&\hspace{1cm}-
 \Lambda_{n+1}^{*} \Sigma_{n+1}^{-1}
( \lambda I_{N}- \alpha )^{-1} \Lambda_{n+1}(I_{2m}-
\frac{i}{\lambda}J).\label{1.5}
\end{eqnarray}
Using the Taylor expansion of $( \lambda I_{N}- \alpha )^{-1}$ at
infinity one shows that (\ref{1.5}) is in its turn equivalent to
the set of equalities:
\begin{eqnarray}
&&S_n-J=\Lambda_n^{*} \Sigma_n^{-1}
 \Lambda_n - \Lambda_{n+1}^{*} \Sigma_{n+1}^{-1}
 \Lambda_{n+1},\label{1.6}\\
\noalign{\vskip6pt} && \Lambda_{n+1}^{*} \Sigma_{n+1}^{-1} \a^p
\Lambda_{n+1}-i \Lambda_{n+1}^{*} \Sigma_{n+1}^{-1} \a^{p-1}
\Lambda_{n+1}J=\nonumber\\
 \noalign{\vskip6pt}
 && \hspace{3cm}= \Lambda_{n}^{*} \Sigma_{n}^{-1} \a^p \Lambda_{n}-i
S_n \Lambda_{n}^{*} \Sigma_{n}^{-1} \a^{p-1} \Lambda_{n} \quad
(p>0).\label{1.7}
\end{eqnarray}
Equality (\ref{1.6}) is equivalent to (\ref{0.6}). Taking into
account the first relation in (\ref{0.7}) we have $\a
\Lambda_{n+1} -i \Lambda_{n+1} J= \a \Lambda_n + \a^{-1}
\Lambda_n$. Thus the equalities in (\ref{1.7}) can be rewritten in
the form $K_n \a^{p-2} \Lambda_n=0$, where
\begin{equation} \label{1.8}
K_n=\Lambda_{n+1}^{*} \Sigma_{n+1}^{-1} (\a^2+I_{N})-
\Lambda_{n}^{*} \Sigma_{n}^{-1} \a^2+i S_n \Lambda_{n}^{*}
\Sigma_{n}^{-1} \a.
\end{equation}
Therefore, if we prove that $K_n=0$, then equalities (\ref{1.7})
will be proved, and so formula (\ref{1.4}) will be proved too.
Substitute (\ref{0.6}) into (\ref{1.8}), and  again use the first
relation in (\ref{0.7}) to obtain
\begin{eqnarray}
K_n&=&\Lambda_{n+1}^{*} \Sigma_{n+1}^{-1} (\a^2+I_{N})-
\Lambda_{n}^{*} \Sigma_{n}^{-1} \a^2+i J \Lambda_{n}^{*}
\Sigma_{n}^{-1} \a +\nonumber\\ \noalign{\vskip6pt} &&  +i
\Lambda_n^* \S_n^{-1} \Lambda_n \Lambda_{n}^{*} \Sigma_{n}^{-1} \a
-i \Lambda_{n+1}^* \S_{n+1}^{-1}(\Lambda_n+i \a^{-1} \Lambda_n
J)\Lambda_{n}^{*} \Sigma_{n}^{-1} \a. \label{1.9}
\end{eqnarray}
Now we use (\ref{1.1}) to obtain $i \Lambda_n \Lambda_{n}^{*}
\Sigma_{n}^{-1}= \a - \S_n \a^* \S_n^{-1}$ and substitute this
relation into (\ref{1.9}). After easy transformations it follows
that
\begin{eqnarray}
K_n&=&\Lambda_{n+1}^{*} \Sigma_{n+1}^{-1} \Big(\a^{-1} \S_n
(\a^*)^{-1}+ \S_n   + \a^{-1} \Lambda_n J \Lambda_n^* (\a^*)^{-1}
\Big) \a^* \S_n^{-1} \a +\nonumber\\ \noalign{\vskip6pt} &&  +i J
\Lambda_{n}^{*} \Sigma_{n}^{-1} \a - \Lambda_n^* \a^* \S_n^{-1}
\a. \label{1.10}
\end{eqnarray}
In view of the second relation in (\ref{0.7}) the first term on
the right hand side of (\ref{1.10}) equals $\Lambda_{n+1}^{*} \a^*
\S_n^{-1} \a$ and we have
\begin{equation} \label{1.11}
K_n=(\Lambda_{n+1}^*+i J \Lambda_n^* (\a^*)^{-1} -\Lambda_n^*)\a^*
\S_n^{-1} \a.
\end{equation}
By the first relation in (\ref{0.7}) the equality $K_n=0$ is now
immediate, i.e., (\ref{1.4}) is true.

Notice  that equality (\ref{1.-1})  is valid for $n=0$. Suppose
that it is valid for  $ n =r$. Then (\ref{0.5}) and (\ref{1.-1})
yield
\begin{equation} \label{1.12}
W_{r+1}(\lambda)=(I_{2m}- \frac{i}{\lambda}S_r)W_{ \alpha ,
\Lambda }(r, \lambda ) \Big( I_{2m}-\frac{i}{\lambda}J \Big)^{r}
W_{ \alpha , \Lambda }(0, \lambda )^{-1}.
\end{equation}
By (\ref{1.4}) and (\ref{1.12}) the validity of (\ref{1.-1}) for
$n=r+1$ easily follows, i.e.,  (\ref{1.-1}) is proved by
induction.

Consider now the matrices $S_n$ given by (\ref{0.6}). It is easy
to see that $S_n=S_n^*$. Notice also that in view of (\ref{1.1})
we have
\begin{equation} \label{1.13}
W_{ \alpha , \Lambda }(r, \lambda )W_{ \alpha , \Lambda }(r,
\overline{ \lambda} )^* =I_{2m} \quad (r \geq 0),
\end{equation}
where $\overline{ \lambda}$ stands for complex conjugate for
$\lambda$. From (\ref{1.4}) and (\ref{1.13}) it follows that
$(I_{2m}- i\lambda^{-1}S_n) (I_{2m}+
i\lambda^{-1}S_n)=\lambda^{-2}(\lambda^2+1)I_{2m}$. Thus the
equality $S_n^*=S_n^{-1}$  holds, which finishes the proof of the
theorem. $\Box$

The case when $\pm i \, \not\in \s(\a)$ ($\s$ means spectrum) is
important for the study of the IHM model. Assume this condition is
fulfilled, and put
\begin{equation}\label{2.21.2}
R_n=(I_N-i \a^{-1})^{-n} \S_n(I_N+i (\a^*)^{-1})^{-n},
\end{equation}
\begin{equation}\label{2.21.2'}
Q_n=(I_N+i \a^{-1})^{-n} \S_n(I_N-i (\a^*)^{-1})^{-n}.
\end{equation}
The following proposition will be useful for formulating the
conditions of invertibility of $\S_n$ in a somewhat different form
then those in Proposition \ref{Pn2.0} (see Corollary \ref{Cy1.2n}
below).
\begin{Pn} \label{Pn1.3}
Let the matrices $\a$ $( \det \, \a \not= 0)$, $\S_0=\S_0^*$, and
$\Lambda_0$ satisfy $($\ref{0.8}$)$, and let the matrices $ \S_n $
be given by $($\ref{0.7}$)$. If $i \not\in \s(\a)$,
 then the
sequence of matrices $\{ R_n \}$ is well defined and
non-decreasing. If $-i \not\in \s(\a)$,
  then the sequence of matrices $\{ Q_n \}$ is well defined
and non-increasing.
\end{Pn}
{\bf Proof.} To prove that the  sequence $\{ R_n \}$  is
non-decreasing it will suffice to show that
\begin{equation} \label{1.14}
\S_{n+1} -(I_N-i \a^{-1}) \S_n(I_N+i (\a^*)^{-1}) \geq 0.
\end{equation}
For this purpose notice that
\begin{eqnarray}
\S_{n+1} &-&(I_N-i \a^{-1}) \S_n(I_N+i (\a^*)^{-1})= \nonumber
\\ &&
=\S_{n+1}- \S_n- \a^{-1} \S_n (\a^*)^{-1}-i \a^{-1}(\a \S_n- \S_n
\a^*) (\a^*)^{-1}. \nonumber
\end{eqnarray}
Hence, in view of (\ref{0.7}) and (\ref{1.1}) we get
\begin{eqnarray}
\S_{n+1} -(I_N-i \a^{-1}) \S_n(I_N&+&i (\a^*)^{-1})= \nonumber \\
&& = \a^{-1}(\Lambda_n J \Lambda_n^* + \Lambda_n \Lambda_n^*)
(\a^*)^{-1}. \label{1.16}
\end{eqnarray}
Since $J+I_{2m} \geq 0$, the inequality (\ref{1.14}) is immediate
from (\ref{1.16}).

Similarly,  from (\ref{0.7}) and (\ref{1.1}) we get
\begin{eqnarray}
\S_{n+1} -(I_N+i \a^{-1}) \S_n(&I_N&-i (\a^*)^{-1})= \nonumber \\
&& = \a^{-1}(\Lambda_n J \Lambda_n^* - \Lambda_n \Lambda_n^*)
(\a^*)^{-1} \leq 0, \label{1.17}
\end{eqnarray}
and so the sequence of matrices $\{ Q_n \}$ is non-increasing.
$\Box$

According to Proposition \ref{Pn1.3}, when  $i \not\in \s(\a)$ and
$\S_0>0$, we have $R_n>0$.
\begin{Cy} \label{Cy1.2n} Let the conditions of Proposition \ref{Pn1.3}
hold, and assume that $i \not\in \s(\a)$ and $\S_0>0$. Then  we
get $\S_n>0$ for all $n>0$.
\end{Cy}
Partition the matrices $W_{ \alpha , \Lambda }(r, \lambda )$ and
$\Lambda_r$ into two $m$-column blocks each:
\[
W_{ \alpha , \Lambda} (r, \lambda )=\Big[ \big( W_{ \alpha ,
\Lambda }(r, \lambda ) \big)_1 \quad \big( W_{ \alpha , \Lambda
}(r, \lambda ) \big)_2 \Big], \quad \Lambda_r= \big[ (\Lambda_r)_1
\quad (\Lambda_r)_2].
\]
The next lemma will be used in Section 4.
\begin{La} \label{La2.3}
Let  the matrices $\a$ $(0, \, \pm i \, \not\in \s(\a))$, $\S_0$
$(\S_0=\S_0^*)$ and $\Lambda_0$ satisfy $($\ref{0.8}$)$, and let
the matrices $ \S_n $ be given by $($\ref{0.7}$)$. Then for $n
\geq 0$ the following relations hold:
\begin{eqnarray}
 \big(& W_{ \alpha , \Lambda }(n,& i ) \big)_1=
 \big( W_{ \alpha , \Lambda }(n+1, -i ) \big)_1
\times \nonumber \\ && \times
  \Big( I_m+2 \big(
W_{ \alpha , \Lambda }(n, i ) \big)_1^* \Lambda_{n}^*
\S_n^{-1}(\a^2+I_N)^{-1}(\Lambda_n)_1 \Big),  \label{1.18}
\end{eqnarray}
\begin{eqnarray}
\big( &W_{ \alpha , \Lambda }(n,& -i ) \big)_2=  \big( W_{ \alpha
, \Lambda }(n+1, i ) \big)_2 \times \nonumber \\ && \times \Big(
I_m-2 \big( W_{ \alpha , \Lambda }(n, -i ) \big)_2^* \Lambda_{n}^*
\S_n^{-1}(\a^2+I_N)^{-1}(\Lambda_n)_2 \Big). \label{1.19}
\end{eqnarray}
\end{La}
{\bf Proof.} From the proof of Theorem \ref{Tm1.1} we know that
$K_n=0$, where $K_n$ is given by (\ref{1.8}). In particular, we
get
\begin{equation} \label{1.20}
K_n \a^{-1}(\a^2+I_N)^{-1}(\Lambda_n)_1=0, \quad K_n
\a^{-1}(\a^2+I_N)^{-1}(\Lambda_n)_2=0.
\end{equation}
To prove (\ref{1.18}) notice that $(\Lambda_{n+1})_1= \a^{-1}(\a+i
I_N) (\Lambda_{n})_1$ and rewrite the first equality in
(\ref{1.20}) as
\begin{eqnarray}
\Lambda_{n+1}^* \S_{n+1}^{-1}(\a &+&i I_N)^{-1}(\Lambda_{n+1})_1-
\Lambda_{n}^* \S_{n}^{-1}(\a -i I_N)^{-1}(\Lambda_{n})_1+
\nonumber \\ && +
 \label{1.21}
i(I_{2m}+S_n) \Lambda_{n}^*
\S_n^{-1}(\a^2+I_N)^{-1}(\Lambda_n)_1=0.
\end{eqnarray}
Put $\lambda = -i$ in (\ref{1.4}) and take into account
(\ref{1.13}) to derive
\begin{equation} \label{1.22}
I_{2m}+S_n=2 \big( W_{ \alpha , \Lambda }(n+1, -i ) \big)_1 \big(
W_{ \alpha , \Lambda }(n, i ) \big)_1^*.
\end{equation}
In view of  definition (\ref{1.0}) of $W_{ \alpha , \Lambda }$,
equality (\ref{1.18}) follows from (\ref{1.21})  and (\ref{1.22}).

Putting $\lambda = i$ in (\ref{1.4}) we get
\begin{equation} \label{1.23}
I_{2m}-S_n=2 \big( W_{ \alpha , \Lambda }(n+1, i ) \big)_2 \big(
W_{ \alpha , \Lambda }(n, -i ) \big)_2^*.
\end{equation}
Analogously to the proof of (\ref{1.18}) we derive (\ref{1.19})
from (\ref{1.23}) and the second equality in (\ref{1.20}). $\Box$

\begin{Rk} \label{Rk2.4}
According to $($\ref{1.22}$)$ and $($\ref{1.23}$)$ the rank of the
matrices $I_{2m} \pm S_n$ is less than or equal to $m$. Together
with the  formula $($\ref{0.5'}$)$ this implies that under the
conditions of Lemma \ref{La2.3} we have $S_n=U_nJU_n^*$, where
$U_n$ are unitary matrices and $J$ is defined in $($\ref{0.1}$)$.
\end{Rk}

\section{Weyl functions: direct and inverse problems}
\label{Weyl} \setcounter{equation}{0} In this section we prove
Theorems \ref{Tm2.1} and \ref{Tm2.3}, and Proposition \ref{Pn2.0}.
At the end of the section a lemma on the case $i \not\in \s(\a)$
is treated too.

\medskip\noindent
{\bf Proof of Proposition \ref{Pn2.0}. }
 Suppose $f$ is an eigenvector of $\a$, that is,
$\a f=c f$, $f \not=0$. Then formula (\ref{2.0}) yields the
equality
\begin{equation} \label{2.0'}
i( \overline{c}-c)f^*f=f^*(\t_1 \t_1^*+\t_2 \t_2^*)f \geq 0.
\end{equation}
So $c \in \overline{\BC_+}$. Moreover, if $c = \overline{c}$, then
according  to (\ref{2.0'}) we have $\t_1^*f= \t_2^*f=0$, and
therefore $\a f= \a^*f$. It follows that
\begin{equation} \label{2.0.2}
f^* \t_1=0,  \quad f^*(\a -cI_N)=0 \quad (f \not=0).
\end{equation}
As $\{\a$, $\t_1\}$ is a full range pair, so the pair $\{\a
-cI_N$, $\t_1 \}$ is full range, which contradicts (\ref{2.0.2}).
This implies  that $c \in \BC_+$, i.e., $\s(\a) \subset \BC_+$.

Recall that identity (\ref{1.1}) was deduced in the proof of
Theorem \ref{Tm1.1}. Taking into account that $\s(\a) \subset
\BC_+$, identity (\ref{1.1}) yields
\begin{equation} \label{2.0.3}
\S_n = \frac{1}{2 \pi} \int_{- \infty}^{\infty} (\a - \lambda
I_N)^{-1} \Lambda_n \Lambda_n^*
 (\a^* - \lambda I_N)^{-1}d \lambda.
\end{equation}
Notice now that the pair $\{\a$, $(I_N+i \a^{-1})^n \t_1 \}$ is
full range and use (\ref{2.14}), (\ref{2.0.3}) to obtain $\S_n>0$
for all $n \geq 0$. $\Box$

\begin{Rk} \label{Rk3.1n}
In the same way as in the proof of Proposition \ref{Pn2.0} above
the inclusion $\s(\a) \subset \BC_+$ follows from the weaker
condition that the pair $\{ \a$, $\Lambda_0 \}$ is full range.
However, the example $N=1$, $\a=i$, $\theta_1=0$, $\t_2 \t_2^*=2$,
which yields $\S_n \equiv 0$ for $n>0$, shows that we have to
require that the pair $\{ \a$, $\t_1 \}$ is full range in order to
get $\S_n>0$. The full range condition on the pair
 $\{\a$, $\Lambda_0 \}$ is not enough for this conclusion.
\end{Rk}

Recall now  Definition \ref{DnAdm} of the admissible triple.
Proposition \ref{Pn2.0} implies, in particular, that $\det \a
\not=0$ for the admissible triple and the spin sequences $\{ S_n
\}$ determined  by it are well defined for all $n \geq 0$. In
other words the class FG  is well defined.

\medskip\noindent
{\bf Proof of Theorem \ref{Tm2.1}.} Let $W_{\alpha , \Lambda}(n,
\lambda)$ be given by (\ref{1.0}).
 Write
$W_{ \alpha , \Lambda }(0, \lambda )$
 as
\begin{equation} \label{2.4}
W_{ \alpha , \Lambda }(0, \lambda )= \left[
\begin{array}{lr}
a( \lambda ) & b( \lambda ) \\ c( \lambda ) & d( \lambda )
\end{array}
\right],
\end{equation}
where
 the $m \times m$ matrix functions $b( \lambda )$ and $d( \lambda
)$ are given by
\begin{equation} \label{2.6}
b( \lambda ) =i \theta_{1}^{*}( \lambda I_{N} - \alpha )^{-1}
\theta_{2}, \hspace{1em}d( \lambda ) =I_{m}+i \theta_{2}^{*}(
\lambda I_{N} - \alpha )^{-1} \theta_{2}.
\end{equation}
We first prove that
\begin{equation} \label{2.5}
b( \lambda )d( \lambda )^{-1} =i \theta_{1}^{*}( \lambda I_{N} -
\beta )^{-1} \theta_{2}.
\end{equation}
Using (\ref{inv}) in preliminaries, from  (\ref{2.3}) and
(\ref{2.6}) we obtain
\begin{equation} \label{2.7}
d( \lambda )^{-1} =I_{m}-i \theta_{2}^{*}( \lambda I_{N} - \beta
)^{-1} \theta_{2}.
\end{equation}
The equalities (\ref{2.3}), (\ref{2.6}) and (\ref{2.7}) yield
\begin{eqnarray}
b( \lambda )d( \lambda )^{-1} = \qquad \qquad \qquad \qquad \qquad
\nonumber \\ =i \theta_{1}^{*}( \lambda I_{N} - \alpha )^{-1}
\theta_{2} +i \theta_{1}^{*}( \lambda I_{N} - \alpha )^{-1}( \beta
- \alpha ) ( \lambda I_{N} - \beta )^{-1} \theta_{2}. \label{2.8}
\end{eqnarray}
From (\ref{2.8}) formula (\ref{2.5}) follows.

Let $\varphi$ be defined by (\ref{2.2}), and thus by virtue of
(\ref{2.5}) we have
\begin{equation} \label{2.9}
\varphi ( \lambda )=b( \lambda )d( \lambda )^{-1}.
\end{equation}
 By (\ref{2.4}), (\ref{2.9})
and the representation (\ref{1.-1}) of the fundamental solution we
get
\begin{equation} \label{2.10}
W_n( \lambda ) \left[ \begin{array}{c} \varphi ( \lambda )
\\ I_{m} \end{array} \right]= \left( \frac{\lambda + i}{\lambda}
\right)^n
 W_{ \alpha , \Lambda }(n, \lambda )
 \left[ \begin{array}{c} 0
\\ d( \lambda )^{-1} \end{array} \right].
\end{equation}
Notice also that (\ref{1.1}) yields a more general formula than
(\ref{1.13}), namely:
\begin{eqnarray}
 &W_{ \alpha , \Lambda }(n,& \lambda )^{*}
 W_{ \alpha , \Lambda }(n, \lambda )= \nonumber \\ &&=
I_{2m}-i( \lambda - \overline{ \lambda })
 \Lambda_n^{*}( \overline{ \lambda }I_{N}- \alpha
^{*})^{-1} \S_n^{-1} (  \lambda I_{N}-  \alpha )^{-1} \Lambda_n.
\label{2.10'}
\end{eqnarray}
As the second term in the right hand side of (\ref{2.10'}) is
nonpositive, it follows  from  formula (\ref{2.10'}) that
\begin{equation} \label{2.11}
 W_{ \alpha , \Lambda }(n, \lambda )^{*}
 W_{ \alpha , \Lambda }(n, \lambda ) \leq
I_{2m} \hspace{1em}( \lambda \in {\BC_{-}}).
\end{equation}
Taking into account (\ref{2.10}) and (\ref{2.11})  we obtain
(\ref{2.1}), i.e., $ \varphi $ is a Weyl function.

It remains only to prove the uniqueness of the Weyl function. Let
us first show that for some $M>0$  and all $n \geq 0$ we have the
inequality
\begin{equation} \label{2.12}
f^*(I_N-i (\a^*)^{-1})^{n} \S_n^{-1}(I_N+i \a^{-1})^{n}f \leq
Mf^*f, \quad f \in L,
\end{equation}
where
\[
L:= {\mathrm {span}}_{\lambda \not\in \s(\a)} {\mathrm{Im}} \,
(\lambda I_N - \a)^{-1} \t_1.
\]
In view of (\ref{2.14}) formula (\ref{2.10'}) yields
\begin{eqnarray}
 \t_1^{*}( \overline{ \lambda }I_{N}- \alpha
^{*})^{-1}
 (&I_N&-i( \a^*)^{-1})^n \S_n^{-1} \times \nonumber \\ && \times
(I_N+i \a^{-1})^n(  \lambda I_{N}-  \alpha )^{-1} \t_1 \leq \frac
{i}{  \overline{ \lambda }- \lambda} I_m. \label{2.15}
\end{eqnarray}
Now to get (\ref{2.12}) from (\ref{2.15}) we note that ${\mathrm
{span}}_{\lambda \not\in \s(\a)} {\mathrm{Im}} \, (\lambda I_N -
\a)^{-1} \t_1$ coincides with the same span when $\lambda$ runs
over an $\varepsilon$-neighbourhood $O_\varepsilon$ of any
$\lambda_0 \not\in \s_{\a}$,  for any sufficiently small
$\varepsilon>0$.

By (\ref{2.10'}) and (\ref{2.12}) we can choose $M_1>0$ such that
 we have
\begin{equation} \label{2.16}
[I_m \quad 0] W_{ \alpha , \Lambda }(n, \lambda )^{*} W_{ \alpha ,
\Lambda }(n, \lambda )
 \left[ \begin{array}{c} I_m
\\ 0 \end{array} \right]    \geq \frac{1}{2} \quad
{\mathrm{for}} \,  {\mathrm{all}} \, | \lambda |>M_1.
\end{equation}
Without loss of generality we may assume that $M_1$ is large
enough in order that $M_1> || \a ||$ and $a(\lambda)$ is
invertible for $| \lambda |>M_1$. Then, taking into account
(\ref{1.-1}) and (\ref{2.16}), we obtain
\begin{eqnarray}
&\sum_{n=0}^{r}&[I_m \qquad \big(c(\lambda)a(\lambda)^{-1}
\big)^*] \times \nonumber \\ && \times W_n(\lambda)^* W_n(\lambda)
\left[
\begin{array}{c}
I_m  \\ c(\lambda)a(\lambda)^{-1}
\end{array}
\right] > \frac{r}{2} \big(a(\lambda)^{-1} \big)^* a(\lambda)^{-1}
\label{2.18}
\end{eqnarray}
for all $\lambda$ in the domain
\[
D= \{ \lambda: \, | \lambda |>M_1, \quad \Im \, \lambda < -
\frac{1}{2} \}.
\]
In other words, for $\lambda \in D$ we have
\begin{equation} \label{2.20}
\sum_{n=0}^{\infty} f^* W_n(\lambda)^* W_n(\lambda)f=  \infty
\quad (f \in L_1),
\end{equation}
where
\[
 L_1:= {\mathrm{Im}} \, \left[
\begin{array}{c}
I_m  \\ c(\lambda)a(\lambda)^{-1}
\end{array}
\right].
\]

Suppose now that $\varphi$ and $ \tilde{ \varphi }$ are Weyl
functions of (\ref{0.5}) and that for some fixed $\lambda_{0} \in
D$
 we have
$ \tilde{ \varphi }( \lambda_{0} ) \not=  \varphi ( \lambda_{0}
)$. Put
\[
L_{2}={ \mathrm{ Im}} \left[ \begin{array}{c} \varphi (
\lambda_{0} )
\\ I_{m} \end{array} \right]
+{ \mathrm{ Im}} \left[ \begin{array}{c} \tilde{ \varphi }(
\lambda_{0} )
\\ I_{m} \end{array} \right].
\]
According to the definition of the Weyl function we have
\begin{equation}\label{2.21}
\sum_{n=0}^{\infty} f^* W_n(\lambda)^* W_n(\lambda)f<  \infty
\quad (f \in L_2).
\end{equation}
As $ \dim L_{1} = m$ and $ \dim L_{2} >m$,  there is a non-zero
vector $f$ such that $f \in (L_{1} \cap L_{2})$, which contradicts
(\ref{2.20}) and (\ref{2.21}). $\Box$

For the proof of Theorem \ref{Tm2.3} we shall use the following
lemma which is of independent interest.

\begin{La}\label{lemreal}
A  strictly proper rational $m \times m$ matrix function $\varphi$
admits a minimal realization of the form
\begin{equation}\label{realphi}
\varphi (\lambda)= i\theta_1^*(\lambda I_N-\beta)^{-1}\theta_2,
\end{equation}
such that
\begin{equation} \label{2.26}
\beta - \beta^{*}=i( \theta_{1} \theta_{1}^{*}- \theta_{2}
\theta_{2}^{*}).
\end{equation}
\end{La}
\medskip\noindent
{\bf Proof.} We may assume that $ \varphi $ is given by the
minimal realization $($\ref{2.24}$)$. First let us show that
equation (\ref{2.23}) has a unique solution $X>0$.

The minimality of the realization (\ref{2.24}) means that the pair
$\{\vt_1^*$, $\g \}$ is observable and the pair $\{ \g$, $\vt_2
\}$ is controllable. Notice that $\im \vt \supseteq \im \vt \vt^*$
and $f^* \vt \vt^*=0$ yields $f^* \vt =0$, i.e., $\im \vt = \im
\vt \vt^*$. Hence the pair $\{ \g$, $\vt_2 \vt_2^* \}$ is
controllable too. Therefore the pair $\{\vt_2 \vt_2^*$, $i \g^*\}$
is observable. The pair $\{i \g^*$, $\vt_1 \}$ is controllable and
hence c-stabilizable, that is, there exists a matrix $K$ such that
$i \g^*+ \vt_1 K$ has all its eigenvalues in the open left
half-plane. But then we can use  Theorem
 16.3.3 \cite{LR} (see also \cite{Kal})
to show that the equation (\ref{2.23}) has a unique non-negative
solution $X$ and that this solution $X$ is positive definite.

Next, let $\t_1$, $\t_2$, $\b$ be defined by (\ref{2.25}). From
(\ref{2.23}) and (\ref{2.25}) it follows that (\ref{2.26}) is
satisfied.

According to (\ref{2.24}) and (\ref{2.25}) the function $\varphi$
is also given by the realization (\ref{realphi}). Moreover as the
realization (\ref{2.24}) is minimal, the same is true for the
realization (\ref{realphi}).   $\Box$

\medskip\noindent
{\bf Proof of Theorem \ref{Tm2.3}.} Let $ \varphi $ be a strictly
proper rational $m \times m$ matrix function. Let  $\t_1$, $\t_2$,
$\b$ be as in Lemma \ref{lemreal}, and put $\a= \b+ i \t_2
\t_2^*$. Then the triple $ \alpha $, $\t_1$, and $\t_2$ satisfies
(\ref{2.0}). Furthermore, the pairs $\{\b$, $\t_2 \}$ and $\{
\b^*$, $\t_1 \}$ are full range. Since $\a=
 \beta +i \theta_{2} \theta_{2}^{*} $, it is
immediate that the pair $\{ \a$, $\t_2 \}$ is full range (see
Section \ref{Ap}). By (\ref{2.0}) we have $\b^*=\a^*+i \t_2
\t_2^*= \a-i \t_1 \t_1^*$. Hence, as  $\{ \b^*$, $\t_1 \}$ is a
full range pair, so the pair $\{\a$, $\t_1 \}$ is full range too.
Therefore the triple $ \alpha $, $ \theta_{1}$, and $\theta_{2}$
is admissible. From Theorem \ref{Tm2.1} it follows now that the
function $\{ S_n \} $ determined by the admissible triple $
\alpha$, $ \theta_{1}$ and $\theta_{2}$ is indeed a solution of
the inverse problem.

Let us prove now the uniqueness of the solution of the inverse
problem. Suppose that there is  system (\ref{0.5}) with another
spin sequence $\{ \wt S_n \}$, given by the admissible triple $\wt
\a$, $\wt \t_1$, $\wt \t_2$, and with the same Weyl function
$\varphi$. According to Theorem \ref{Tm2.1} we have another
realization for $\varphi$, namely
\begin{equation} \label{2.31}
\varphi ( \lambda ) =i \wt \theta^{*}_{1}( \lambda I_{\wt N}- \wt
\b)^{-1}\wt \theta_{2}, \quad \wt \b= \wt \a- i \wt \t_2 \wt
\t_2^*.
\end{equation}

As the pairs $\{ \wt \a$, $\wt \t_1 \}$ and $\{ \wt \a$, $\wt \t_2
\}$ are controllable, and $\wt \b=\wt \a - i \wt \t_2 \wt \t_2^*$,
$\wt \b^*= \wt a-i \wt \t_1 \wt \t_1^*$, it follows that the pairs
$\{ \wt \b$, $\wt \t_2 \}$ and  $\{ \wt \b^*$, $\wt \t_1 \}$ are
also controllable. Thus the realization (\ref{2.31}) is minimal
and $\wt N =N$. Moreover, there is (see Section \ref{Ap}) a state
space similarity transforming the realization (\ref{2.2}) into
(\ref{2.31}), that is, there exists an invertible $N \times N$
matrix $S$ such that
\begin{equation} \label{2.32}
\wt \b =S \b S^{-1}, \quad \wt \t_2 = S \t_2, \quad \wt \t_1^*=
\t_1^* S^{-1}.
\end{equation}
Recall that $\wt \a$, $\wt \t_1$, $\wt \t_2$ is an admissible
triple, and therefore we have
\begin{equation} \label{2.33}
\wt \beta - \wt \beta^{*}=i(\wt \theta_{1} \wt \theta_{1}^{*}- \wt
\theta_{2} \wt \theta_{2}^{*}).
\end{equation}
From (\ref{2.32}) and (\ref{2.33}) it follows that
$Z=S^{-1}(S^*)^{-1}$ satisfies
\begin{equation} \label{2.29}
\b Z - Z \b^*=i(Z \t_1 \t_1^* Z- \t_2 \t_2^*), \quad Z>0.
\end{equation}
Completely similar to the uniqueness of $X>0$ in (\ref{2.23}) one
obtains the uniqueness of the solution $Z>0$ of (\ref{2.29}).
Compare now (\ref{2.26}) and (\ref{2.29}) to see that $Z=I_N$ and
thus $S$ is unitary. In view of
 this we have
\[
\wt \a =S \a S^*, \quad \wt \t_2 = S \t_2, \quad \wt \t_1= S \t_1.
\]
The last transformation does not change the spin sequence, i.e.,
$\{ S_n \}= \{ \wt S_n \}$. $\Box$

The case when $i \not\in \s(\a)$ is important in the next section.
We shall use the acronym $\wt{\mathrm{FG}}$ to denote the class of
spin sequences $\{ S_n\}$ determined by the triples $\a$,  $\t_1$,
$\t_2$, with $\a$ an $N \times N$ non-singular matrix and $\t_1$,
$\t_2$ of size $N \times m$, satisfying the identity (\ref{2.0})
and the additional special condition $ i \not\in \s(\a)$. Notice
(see the beginning of the proof of Proposition \ref{Pn2.0}) that
(\ref{2.0}) implies that $\s(\a) \subset \overline{\BC_+}$. The
next lemma shows that without loss of generality we can also
require that $\{ \a$, $\t_1\}$ and $\{ \a$, $\t_2 \}$ are full
range pairs, i.e., $\wt{\mathrm{FG}} \, \subseteq$ FG.
\begin{La} Assume that the spin sequence
$\{ S_n \}$ belongs $\wt{\mathrm{FG}}$. Then it can be determined
by a triple $\a$, $\t_1$, $\t_2$ such that $\a$ is non-singular,
$($\ref{2.0}$)$ holds,
 $i \, \not\in \s(\a)$, and the pairs
$\{\a$, $\t_1 \}$ and $\{\a$, $\t_2 \}$ are full range.
\end{La}
{\bf Proof.} Let $N$ denote the minimal order of $\a$ ($0, \,  i
\, \not\in \s(\a)$) in the set of triples that satisfy (\ref{2.0})
and determine the given spin sequence $\{S_n \}$. Suppose the $N
\times N$ matrix $\wh \a$ and the $N \times m$ matrices $\wh \t_1$
and $\wh \t_2$ form such a triple but the pair $\{ \wh \a$, $\wh
\t_2 \}$ is not full range. Put
\[ {\displaystyle{\wh L_0: =  {\mathrm {span}} \,
\bigcup_{k=0}^{\infty} {\mathrm{Im}} \, \wh \a^k \wh \t_2, \,
N_0:= \dim \widehat L_0}} \] and choose a unitary matrix $q$ that
maps $\wh L_0$ onto the $ L_0:= {\mathrm{Im}} \, [I_{N_0} \quad
0]^*$. Then we have
\begin{equation}\label{2.21.1}
 \a:=q \wh \a q^*=
 \left[ \begin{array}{lr}\wt  \a &  \a_{12} \\
0 &  \a_{22} \end{array} \right], \quad  \t_1:= q \wh \t_1= \left[
\begin{array}{c}\wt  \t_{1}
\\ \kappa \end{array} \right], \quad
 \t_2:= q \wh \t_2=
 \left[ \begin{array}{c} \wt \t_{2}
\\ 0 \end{array} \right],
\end{equation}
where the  $N_0 \times N_0$ matrix $\wt \a$ ($0, \,  i \, \not\in
\s(\wt \a)$) and the $N_0 \times m$ matrices $\wt \t_{1}$,
 $\wt
\t_{2}$ form a triple which satisfies  (\ref{2.0}) and determines
$\{S_n \}$. To show this we  need to make some preparations. As
$\wh L_0$ is an invariant subspace of $\wh \a$, so $ L_0$ is an
invariant subspace of $ \a$, and thus $ \a$ has the block
triangular form given in (\ref{2.21.1}). Moreover in view of the
inclusion Im $\wh \t_2 \subseteq \wh L_0$, we have Im $ \t_2
\subseteq  L_0$, i.e.,
 $ \t_2$ has the block form given in
(\ref{2.21.1}). Taking into account that $q$ is unitary, $0, \,  i
\, \not\in \s(\wh \a)$
 and that
$\wh \a$, $\wh \t_1$, $\wh \t_2$  satisfy (\ref{2.0}) and
determine $\{ S_n \}$, we see that $0, \,  i \, \not\in \s( \a)$
and that
 $ \a$, $ \t_1$, $ \t_2$ satisfy
(\ref{2.0})  and determine $\{ S_n \}$ too. So in view of
(\ref{2.21.1}) we have $0, \,  i \, \not\in \s(\wt \a)$ and the
triple $\wt \a$,
 $\wt \t_{1}$,  $\wt \t_{2}$  satisfies (\ref{2.0})
as well.
 Now use the matrices $Q_n$ defined in (\ref{2.21.2'}) to
rewrite $\Lambda_n^* \S_n^{-1} \Lambda_n$ as a $2\times 2$ block
matrix with block of size $m\times m$. This yields
\[
\Lambda_n^* \S_n^{-1} \Lambda_n= \{ (\Lambda_n^* \S_n^{-1}
\Lambda_n)_{kj} \}_{k,j=1}^2
\]
with blocks
\begin{eqnarray} \nonumber
(\Lambda_n^* \S_n^{-1} \Lambda_n)_{11}&=& \t_1^* Q_n^{-1} \t_1 ,
\\ \nonumber
(\Lambda_n^* \S_n^{-1} \Lambda_n)_{12}&=& \t_1^*Q_n^{-1}(I_N+i
\a^{-1})^{-n} (I_N-i \a^{-1})^{n} \t_2,
\\\label{2.21.3}
(\Lambda_n^* \S_n^{-1} \Lambda_n)_{21}&=&  \t_2^*(I_N+i
(\a^*)^{-1})^{n}(I_N-i (\a^*)^{-1})^{-n}Q_n^{-1} \t_1,
\\
(\Lambda_n^* \S_n^{-1} \Lambda_n)_{22}&=& \t_2^*(I_N+i
(\a^*)^{-1})^{n}(I_N-i (\a^*)^{-1})^{-n}Q_n^{-1}\times \nonumber
\\
&&\hspace{2.5cm} \times (I_N+i \a^{-1})^{-n} (I_N-i \a^{-1})^{n}
\t_2. \nonumber
\end{eqnarray}
Partition  $Q_{n+1}-Q_n$ into four blocks:
$Q_{n+1}-Q_n=\{\chi_{kj} \}_{k,j=1}^2$, where $\chi_{11}$ is an
$N_0 \times N_0$ block.  In view of (\ref{2.14}), (\ref{1.17}),
and (\ref{2.21.1}) we obtain
\begin{eqnarray} \chi_{21}&=&0, \quad
\chi_{12}=0, \quad \chi_{22}=0, \quad \chi_{11}= -2 \wt
\a^{-1}(I_{N_0}+i \wt \a^{-1})^{-n-1} \times \nonumber
\\ &&
\times (I_{N_0}-i \wt \a^{-1})^{n} \wt \t_2 \wt \t_2^*(I_{N_0}+i
(\wt \a^*)^{-1})^{n}(I_{N_0}-i (\wt \a^*)^{-1})^{-n-1} (\wt
\a^*)^{-1} . \nonumber
\end{eqnarray}
Denote by $\wt \Lambda_n$, $\wt \S_n$, $\wt Q_n$, $\wt S_n$ et
cetera the matrices generated by the triple  $\wt \a$,
 $\wt \t_{1}$,  $\wt \t_{2}$. One can see that
 $\chi_{11}=\wt Q_{n+1}- \wt Q_n$.
Taking into account $Q_0=I_N$ we obtain that the matrices $Q_n$
are block diagonal: $Q_n = {\mathrm{diag}} \, \{ \wt Q_n, \,
I_{N-N_0} \}$. Now according to (\ref{2.21.1}), (\ref{2.21.3}) it
follows that
\begin{equation}\label{2.21.5}
\Lambda_n^* \S_n^{-1} \Lambda_n- \wt \Lambda_n^* \wt \S_n^{-1} \wt
\Lambda_n= \left[
\begin{array}{lr}
\kappa^* \kappa & 0 \\ 0 & 0
\end{array}
\right].
\end{equation}
By (\ref{0.6}) and (\ref{2.21.5}) we get $S_n= \wt S_n$, i.e., the
triple  $\wt \a$,
 $\wt \t_{1}$,  $\wt \t_{2}$ determines the spin
sequence $\{ S_n \}$. This contradicts the assumption that $N$ is
minimal and therefore the pair $\{\wh \a$, $\wh \t_2 \}$ should be
full range.

In the same way we shall show that the pair $\{\wh \a$, $\wh \t_1
\}$ is full range too. Indeed, suppose $\{\wh \a$, $\wh \t_1 \}$
is not full range. Put now
\[ {\displaystyle{\wh L_0: =  {\mathrm {span}} \,
\bigcup_{k=0}^{\infty} {\mathrm{Im}} \, \wh \a^k \wh \t_1}}, \quad
{\displaystyle{N_0:= \dim \wh L_0}} \]
 and choose a unitary matrix $q$ that maps $\wh
L_0$ onto the $ L_0:= {\mathrm{Im}} \, [I_{N_0} \quad 0]^*$. Then
similar to the previous case we obtain
\begin{equation}\label{2.21.6}
 \a:=q \wh \a q^*=
 \left[ \begin{array}{lr}\wt  \a &  \a_{12} \\
0 &  \a_{22} \end{array} \right], \quad  \t_1:= q \wh \t_1= \left[
\begin{array}{c}\wt  \t_{1}
\\ 0 \end{array} \right], \quad
 \t_2:= q \wh \t_2=
 \left[ \begin{array}{c} \wt \t_{2}
\\ \kappa \end{array} \right],
\end{equation}
where the $N_0 \times N_0$ matrix $\wt \a$ ($0, \,  i \, \not\in
\s(\wt \a)$) and the $N_0 \times m$ matrices $\wt \t_{1}$,
 $\wt
\t_{2}$ form a triple which satisfies (\ref{2.0}).  To show that
the triple $\wt \a$, $\wt \t_{1}$, $\wt \t_{2}$ determines $S_n$
 we shall use the fact that $ \a$, $ \t_{1}$,  $
\t_{2}$ determines $S_n$ and rewrite $\Lambda_n^* \S_n^{-1}
\Lambda_n$ as a $2\times 2$ block matrix with blocks of size
$m\times m$ as follows:
\[
\Lambda_n^* \S_n^{-1} \Lambda_n= \{ (\Lambda_n^* \S_n^{-1}
\Lambda_n)_{kj} \}_{k,j=1}^2
\]
with
\begin{eqnarray} \nonumber
(\Lambda_n^* \S_n^{-1} \Lambda_n)_{22}&=& \t_2^* R_n^{-1} \t_2 ,
\\  \nonumber
 (\Lambda_n^* \S_n^{-1} \Lambda_n)_{21}&=& \t_2^*R_n^{-1}(I_N-i
\a^{-1})^{-n} (I_N+i \a^{-1})^{n} \t_1,
\\
\label{2.21.7} (\Lambda_n^* \S_n^{-1} \Lambda_n)_{12}&=&
\t_1^*(I_N-i (\a^*)^{-1})^{n}(I_N+i (\a^*)^{-1})^{-n}R_n^{-1}
\t_2, \\ (\Lambda_n^* \S_n^{-1} \Lambda_n)_{11}&=& \t_1^*(I_N-i
(\a^*)^{-1})^{n}(I_N+i (\a^*)^{-1})^{-n}R_n^{-1}\times \nonumber
\\
\nonumber &&\hspace{2.5cm} \times (I_N-i \a^{-1})^{-n} (I_N+i
\a^{-1})^{n} \t_1.
\end{eqnarray}
Here the matrices $R_n$ are given by (\ref{2.21.2}). Partition now
$R_{n+1}-R_n$ into four blocks: $R_{n+1}-R_n=\{\chi_{kj}
\}_{k,j=1}^2$, where $\chi_{11}$ is an $N_0 \times N_0$ block.  In
view of (\ref{1.16}) and (\ref{2.21.6}) we obtain
\begin{eqnarray} \nonumber
\chi_{21}&=&0, \quad \chi_{12}=0, \quad \chi_{22}=0, \quad
\chi_{11}= 2 \wt \a^{-1}(I_{N_0}-i \wt \a^{-1})^{-n-1} \times
\\ &&
\times (I_{N_0}+i \wt \a^{-1})^{n} \wt \t_1 \wt \t_1^*(I_{N_0}-i
(\wt \a^*)^{-1})^{n}(I_{N_0}+i (\wt \a^*)^{-1})^{-n-1} (\wt
\a^*)^{-1} . \nonumber
\end{eqnarray}
Denote by $\wt \Lambda_n$, $\wt \S_n$, $\wt R_n$, $\wt S_n$ et
cetera the matrices generated by the triple  $\wt \a$,
 $\wt \t_{1}$,  $\wt \t_{2}$. One can see that
$\chi_{11}=\wt R_{n+1}-
 \wt R_n$.
Taking into account $R_0=I_N$ we obtain that the matrices $R_n$
are block diagonal: $R_n = {\mathrm{diag}} \, \{ \wt R_n, \,
I_{N-N_0} \}$. Now according to (\ref{2.21.6}), (\ref{2.21.7}) it
follows that
\begin{equation}\label{2.21.9}
\Lambda_n^* \S_n^{-1} \Lambda_n- \wt \Lambda_n^* \wt \S_n^{-1} \wt
\Lambda_n= \left[
\begin{array}{lr} 0 & 0 \\ 0&
\kappa^* \kappa
\end{array}
\right].
\end{equation}
By (\ref{0.6}) and (\ref{2.21.9}) we get $S_n= \wt S_n$, i.e., the
triple  $\wt \a$,
 $\wt \t_{1}$,  $\wt \t_{2}$ determines the spin
sequence $\{ S_n \}$. So the pair $\{\wh \a$, $\wh \t_1 \}$ should
also be full range. $\Box$

Finally,  from the proof of Theorem \ref{Tm2.1} we have the
following corollary.
\begin{Cy} \label{Cy3.4}
Let the parameter matrices $\a$ $(0, \,  i \, \not\in \s(\a))$,
$\S_0(0)>0$ and $\Lambda_0(0)$ satisfy the identity
$($\ref{0.8}$)$. Then the Weyl function $\varphi$ of the system
determined by these matrices is given by the formula
\[
\varphi (\lambda)=i \t_1^* \S_0^{-1}(\lambda I_N - \wt \b)^{-1}
\t_2, \quad \wt \b =\a -i \t_2 \t_2^* \S_0^{-1}.
\]
\end{Cy}
This corollary is proved by transforming the matrices $\a$,
$\S_0(0)$ and $\Lambda_0(0)$ into the equivalent set $\S_0^{-
\frac{1}{2}} \a \S_0^{ \frac{1}{2}}$, $I_N$, $\S_0^{- \frac{1}{2}}
\Lambda_0$.

\section{Isotropic Heisenberg magnet} \label{IHM}
\setcounter{equation}{0}

Explicit solutions of the  discrete integrable nonlinear equations
form an interesting and actively studied domain (see  references
in \cite{AS,EJ,FT,KN,KSS,LRag,MS,T}). To study the IHM model we
insert an additional variable $t$ in our notations:
$\Lambda_n(t)$, $\S_n(t)$, $S_n(t)$, $W_{ \alpha , \Lambda }(n,t,
\lambda )$, $\varphi (t, \lambda)$ and so on. Notice that the
order $N$ and the parameter matrix $\a$ do not depend on $t$, and
the dependence on $t$ of the other matrix functions is defined by
the equations
\begin{eqnarray} \label{4.1}
\frac{d \Lambda_0}{dt}&=&-2 \big( (\a -iI_N)^{-1} \Lambda_0 P_+ +
(\a +iI_N)^{-1} \Lambda_0 P_- \big), \\ P_{ \pm
}&=&\frac{1}{2}(I_{2m} \pm J), \nonumber
\end{eqnarray}
\begin{eqnarray}
\frac{d \S_0}{dt}&=&- \Big( (\a -iI_N)^{-1} \S_0(t)+(\a
+iI_N)^{-1} \S_0(t)+ \nonumber
\\  \label{4.2}
 && + \S_0(t)(\a^* +iI_N)^{-1}+ \S_0(t)(\a^* -iI_N)^{-1}+
\\
\nonumber &&+ 2(\a^2+I_N)^{-1} \big( \a \Lambda_0(t)J
\Lambda_0(t)^*+\Lambda_0(t)J \Lambda_0(t)^* \a^* \big)(
(\a^*)^2+I_N)^{-1} \Big).
\end{eqnarray}
We assume that the parameter matrices $\a$, $\S_0(0)$ and
$\Lambda_0(0)$ satisfy the identity
\begin{equation} \label{4.3}
\a \S_0(t)- \S_0(t) \a^*= i \Lambda_0(t) \Lambda_0(t)^*
\end{equation}
at $t=0$. Then according to (\ref{4.1}), (\ref{4.2}) the identity
(\ref{4.3}) holds for all $t$. (This result can be obtained by
differentiating both sides of (\ref{4.3}).)
\begin{Tm} \label{Tm4.1}
Assume the parameter matrices $\a$ $(0, \,  i \, \not\in \s(\a))$,
$\S_0(0)>0$ and $\Lambda_0(0)$ satisfy the identity
\[
\a \S_0(0)- \S_0(0) \a^*= i \Lambda_0(0) \Lambda_0(0)^*.
\]
Define $\S_0(t)$ and $\Lambda_0(t)$ by  equations $($\ref{4.1}$)$
and $($\ref{4.2}$)$. Then $\S_0(t)>0$ on some interval $-
\varepsilon <t < \varepsilon$, and the sequence $\{ S_n(t) \}$
given by $($\ref{0.6}$)$ and $($\ref{0.7}$)$ belongs to
$\widetilde {\mathrm{FG}}$ for each $t$ from this interval.
Moreover, $\{ S_n(t) \}$ $(- \varepsilon <t < \varepsilon )$
satisfies the IHM equations $($\ref{0.5.2}$)$, $($\ref{0.5.1}$)$.
\end{Tm}
{\bf Proof.} As $\a$ does not depend on $t$ and (\ref{4.3}) is
true, it is immediate that $\{ S_n(t) \} \subset
\widetilde{\mathrm{FG}}$.

Similar to the cases treated in \cite{SaA2,SaA3} we shall
successively obtain the derivatives $\frac{d}{d t} \Lambda_n$,
$\frac{d}{d t} \S_n$, $\frac{d}{d t} \big( \Lambda_n^* \S_n^{-1}
\big)$, and $\frac{d}{d t} W_{\a,\Lambda}(n, t, \lambda)$, and use
the expressions for these derivatives to derive the zero curvature
equation (\ref{0.5.-2}), which is equivalent to (\ref{0.5.2}),
(\ref{0.5.1}); see \cite{FT}. In view of (\ref{2.14}) and
(\ref{4.1}) we have
\begin{equation} \label{4.5}
\Lambda_n(t)= \big[(I_N+i \a^{-1})^n e^{-2t(\a - iI_N)^{-1}} \t_1
\quad (I_N-i \a^{-1})^n e^{-2t(\a + iI_N)^{-1}} \t_2 \big],
\end{equation}
where $\t_p:=(\Lambda_0(0))_p$ ($p=1,2$). Hence it follows that
\begin{equation} \label{4.6}
\frac{d \Lambda_n}{dt}=-2 \big( (\a -iI_N)^{-1} \Lambda_n P_+ +
(\a +iI_N)^{-1} \Lambda_n P_- \big).
\end{equation}
Now we shall show by induction that
\begin{eqnarray}
\frac{d \S_n}{dt}&=&- \Big( (\a -iI_N)^{-1} \S_n(t)+(\a
+iI_N)^{-1} \S_n(t)+ \nonumber
\\
\label{4.7} && + \S_n(t)(\a^* +iI_N)^{-1}+ \S_n(t)(\a^*
-iI_N)^{-1}+
\\
\nonumber && + 2(\a^2+I_N)^{-1} \big( \a \Lambda_n(t)J
\Lambda_n(t)^*+\Lambda_n(t)J \Lambda_n(t)^* \a^* \big)(
(\a^*)^2+I_N)^{-1} \Big).
\end{eqnarray}
By (\ref{4.2}) formula (\ref{4.7}) is true for $n=0$. Suppose it
is true for $n=r$. Then, taking into account the second relation
in (\ref{0.7}) we obtain
\begin{eqnarray}
\frac{d \S_{r+1}}{dt}&=&- \Big( (\a -iI_N)^{-1} \S_{r+1}(t) -(\a
-iI_N)^{-1} \a^{-1} \Lambda_r(t)J \Lambda_r(t)^* (\a^*)^{-1}+
\nonumber
\\ \nonumber && +
(\a +iI_N)^{-1} \S_{r+1}(t)-(\a +iI_N)^{-1} \a^{-1} \Lambda_r(t)J
\Lambda_r(t)^* (\a^*)^{-1}+ \\ \nonumber && +
 \S_{r+1}(t)(\a^* +iI_N)^{-1}-\a^{-1} \Lambda_r(t)J \Lambda_r(t)^*
(\a^*)^{-1}(\a^* +iI_N)^{-1}+ \\ \nonumber && + \S_{r+1}(t)(\a^*
-iI_N)^{-1}- \a^{-1} \Lambda_r(t)J \Lambda_r(t)^* (\a^*)^{-1}(\a^*
-iI_N)^{-1}+ \\  &&+ 2(\a^2+I_N)^{-1} \big( (\a+ \a^{-1})
\Lambda_r(t)J \Lambda_r(t)^*+\Lambda_r(t)J \Lambda_r(t)^* \times
\label{4.8}
\\ \nonumber && \times
 (\a^* + (\a^*)^{-1}) \big)
 ( (\a^*)^2+I_N)^{-1} \Big)+ \frac{d}{dt}\big(\a^{-1}
\Lambda_r(t)J \Lambda_r(t)^* (\a^*)^{-1} \big) .
\end{eqnarray}
In view of (\ref{4.6}) we easily calculate that
\begin{eqnarray}\label{4.9}
C_1(t)&:&=     \frac{d}{dt}\big(\a^{-1} \Lambda_r(t)J
\Lambda_r(t)^* (\a^*)^{-1} \big) +
\\ && +(\a -iI_N)^{-1} \a^{-1} \Lambda_r(t)J
\Lambda_r(t)^*
 (\a^*)^{-1}+ (\a +iI_N)^{-1} \a^{-1} \Lambda_r(t)J \times
\nonumber \\ \nonumber && \times  \Lambda_r(t)^* (\a^*)^{-1}
+\a^{-1} \Lambda_r(t)J \Lambda_r(t)^* (\a^*)^{-1}(\a^*
+iI_N)^{-1}+ \a^{-1} \times \\ \nonumber && \times
 \Lambda_r(t)J \Lambda_r(t)^* (\a^*)^{-1}(\a^* -iI_N)^{-1}=  -(\a
-iI_N)^{-1} \a^{-1} \Lambda_r(t)  \times \\ \nonumber && \times
\Lambda_r(t)^* (\a^*)^{-1}+(\a +iI_N)^{-1} \a^{-1} \Lambda_r(t)
\Lambda_r(t)^* (\a^*)^{-1}- \a^{-1} \Lambda_r(t) \times
\\ \nonumber && \times
 \Lambda_r(t)^* (\a^*)^{-1}(\a^* +iI_N)^{-1}+ \a^{-1}
\Lambda_r(t) \Lambda_r(t)^* (\a^*)^{-1}(\a^* -iI_N)^{-1}.
\end{eqnarray}
Notice that  $(\a -iI_N)^{-1}-(\a +iI_N)^{-1}=2i (\a^2
+I_N)^{-1}$. Therefore we can rewrite (\ref{4.9}) as
\begin{eqnarray}
C_1(t)&=&2i \big(  \a^{-1} \Lambda_r(t) \Lambda_r(t)^*
(\a^*)^{-1}((\a^*)^2 +I_N)^{-1}- \nonumber \\  \label{4.10} &&
 -(\a^2 +I_N)^{-1}\a^{-1} \Lambda_r(t) \Lambda_r(t)^* (\a^*)^{-1}
 \big)
.
\end{eqnarray}
Notice also that
\begin{eqnarray}
C_2(t)&:&=2(\a^2+I_N)^{-1} \big( (\a+ \a^{-1}) \Lambda_r(t)J
\Lambda_r(t)^*+ \Lambda_r(t)J \Lambda_r(t)^* \times \nonumber \\
\nonumber && \times (\a^* + (\a^*)^{-1}) \big) (
(\a^*)^2+I_N)^{-1} =2 \big( \a^{-1} \Lambda_r(t)J \Lambda_r(t)^*
\times
\\ && \times ( (\a^*)^2+I_N)^{-1}+ (\a^2+I_N)^{-1} \Lambda_r(t)J
\Lambda_r(t)^* (\a^*)^{-1} \big). \label{4.11}
\end{eqnarray}
Using the first relation in (\ref{0.7}) and equalities
(\ref{4.10}), (\ref{4.11}) we get
\begin{eqnarray} \nonumber
C_2(t)-C_1(t)=2 \big( \a^{-1} \Lambda_r(t)J \Lambda_{r+1}(t)^* (
(\a^*)^2+I_N)^{-1}+ \qquad
\\  \nonumber \,
+(\a^2+I_N)^{-1}  \Lambda_{r+1}(t)J \Lambda_r(t)^* (\a^*)^{-1}
\big)= 2(\a^2+I_N)^{-1} \times
\\
 \label{4.12} \,
\times \big( \a \Lambda_{r+1}(t)J
\Lambda_{r+1}(t)^*+\Lambda_{r+1}(t)J \Lambda_{r+1}(t)^* \a^*
\big)( (\a^*)^2+I_N)^{-1}.
\end{eqnarray}
From (\ref{4.8}) and (\ref{4.12}) it follows that (\ref{4.7}) is
valid for $n=r+1$ and so for all $n>0$.

Taking into account (\ref{4.6}) and (\ref{4.7}) we can obtain the
equation
\begin{eqnarray} \nonumber
\frac{d}{dt} \big( \Lambda_n  (t)^* \S_n(t)^{-1} \big)&=& H_n^+(t)
\Lambda_n
  (t)^* \S_n(t)^{-1}(\a - i I_N)^{-1}+
\\ &&
 \label{4.13}
+H_n^-(t)\Lambda_n(t)^* \S_n(t)^{-1}(\a +i I_N)^{-1},
\end{eqnarray}
where
\begin{equation} \label{4.14}
H_n^+(t)=2W_{\a, \Lambda}(n,t,i)P_+W_{\a, \Lambda}(n,t,-i)^*,
\end{equation}
\begin{equation} \label{4.15}
H_n^-(t)=2W_{\a, \Lambda}(n,t,-i)P_-W_{\a, \Lambda}(n,t,i)^*.
\end{equation}
Indeed, by (\ref{4.6}) we have
\begin{eqnarray} \nonumber
\frac{d}{dt} \big(& \Lambda_n & (t)^* \S_n(t)^{-1} \big)= - \Big(
2P_+ \Lambda_n(t)^*(\a^*+i I_N)^{-1}+ \\ && + 2P_-
\Lambda_n(t)^*(\a^*-i I_N)^{-1} +
 \label{4.16}
 \S_n(t)^{-1} \frac{ d \S_n}{dt}(t) \Big) \S_n(t)^{-1}.
\end{eqnarray}
Identity (\ref{1.1}) yields
\begin{eqnarray} \nonumber
((\a^*) \pm i I_N)^{-1}\S_n(t)^{-1}=\S_n(t)^{-1}(\a \pm i
I_N)^{-1}+ \qquad
\\ \quad
 \label{4.17}
+i ((\a^*) \pm i I_N)^{-1}\S_n(t)^{-1} \Lambda_n(t) \Lambda_n(t)^*
\S_n(t)^{-1}(\a \pm i I_N)^{-1}.
\end{eqnarray}
Finally notice that
\begin{equation} \label{4.18}
2( \a^2 + I_N)^{-1} \a= (\a - i I_N)^{-1}+ (\a +i I_N)^{-1}.
\end{equation}
By using (\ref{4.7}), (\ref{4.17}) and (\ref{4.18}), and after
some calculations, we rewrite (\ref{4.16}) in the form
(\ref{4.13}), where
\begin{eqnarray} \nonumber H_n^+(t)&=& I_{2m}+J-i
\Lambda_n(t)^*\S_n(t)^{-1}(\a - i I_N)^{-1}\Lambda_n(t)J+ \\
\nonumber && + i J \Lambda_n(t)^*(\a^* - i I_N)^{-1}
\S_n(t)^{-1}\Lambda_n(t) + \Lambda_n(t)^*\S_n(t)^{-1} \times
\\ \label{4.19} && \times (\a - i I_N)^{-1}\Lambda_n(t)J \Lambda_n(t)^*(\a^*
- i I_N)^{-1}\S_n(t)^{-1}\Lambda_n(t),
\end{eqnarray}
\begin{eqnarray} \nonumber
H_n^-(t)&=&I_{2m}-J+i \Lambda_n(t)^*\S_n(t)^{-1}(\a + i
I_N)^{-1}\Lambda_n(t)J-
\\ \nonumber && -
i J \Lambda_n(t)^*(\a^* + i I_N)^{-1}\S_n(t)^{-1}\Lambda_n(t) -
\Lambda_n(t)^*\S_n(t)^{-1} \times \\ \label{4.20} && \times (\a +
i I_N)^{-1}\Lambda_n(t)J \Lambda_n(t)^*(\a^* + i
I_N)^{-1}\S_n(t)^{-1}\Lambda_n(t).
\end{eqnarray}
From (\ref{1.0}) and (\ref{4.19}) it follows that
\[ H_n^+(t)=I_{2m}+W_{\a, \Lambda}(n,t,i)JW_{\a, \Lambda}(n,t,-i)^*, \]
 and so,
taking into account (\ref{1.13}), we derive (\ref{4.14}). Equality
(\ref{4.15}) follows from (\ref{4.20}) in a similar way.

Recall now that $m=2$ and (\ref{1.13}) holds. Then according to
(\ref{1.18}), (\ref{1.22}) and (\ref{4.14}) we get
\begin{equation} \label{4.21}
H_n^+(t)=c_n^+(t)(I_2+S_n(t))(I_2+S_{n-1}(t)),
\end{equation}
and according to (\ref{1.19}), (\ref{1.23}) and (\ref{4.15}) we
get
\begin{equation} \label{4.22}
H_n^-(t)=c_n^-(t)(I_2-S_n(t))(I_2-S_{n-1}(t)),
\end{equation}
where $c_n^{\pm}(t)$ are scalar functions. In view of (\ref{4.14})
and (\ref{4.15}) we obtain also Tr $H_n^{\pm}$, where Tr denotes
the trace. Indeed,
\begin{equation} \label{4.23}
{\mathrm{Tr}} \, H_n^{\pm}(t) \equiv 2.
\end{equation}
By Remark \ref{Rk2.4} and formula (\ref{0.5.2}) we derive
\begin{eqnarray} \nonumber
{\mathrm{Tr}} \, (I_2 \pm S_n(t))(I_2 \pm
S_{n-1}(t)&)&={\mathrm{Tr}} \, (I_2 + S_n(t) S_{n-1}(t)) = \\ &&
=2 (1+ \overrightarrow{S}_{n-1}(t) \, \cdot \,
\overrightarrow{S}_{n}(t)).  \label{4.24}
\end{eqnarray}
From (\ref{4.21})-(\ref{4.24}) it follows that $1+
\overrightarrow{S}_{n-1}(t) \, \cdot \, \overrightarrow{S}_{n}(t)
\not= 0$. Now formulas (\ref{0.5.4}), (\ref{4.24}) and
(\ref{4.23}) yield that
\begin{equation} \label{4.25}
{\mathrm{Tr}} \, V_n^{\pm}(t) \equiv 2 \equiv {\mathrm{Tr}} \,
H_n^{\pm}(t).
\end{equation}
Taking into account (\ref{0.5.4}), (\ref{4.21}) and (\ref{4.22})
we see that $V_n^{\pm}= \widehat c^{\pm} H_n^{\pm}$, and so
(\ref{4.25}) yields equalities $V_n^{\pm} \equiv H_n^{\pm}$. Hence
we have
\begin{eqnarray} \nonumber
\frac{d}{dt} \big( \Lambda_n(t)^* \S_n(t)^{-1} \big)&=&
V_n^+(t)\Lambda_n(t)^* \S_n(t)^{-1}(\a - i I_N)^{-1}+
\\
 \label{4.26} &&+
V_n^-(t)\Lambda_n(t)^* \S_n(t)^{-1}(\a +i I_N)^{-1}.
\end{eqnarray}
From $V_n^{\pm} \equiv H_n^{\pm}$, (\ref{1.13}) and definitions
(\ref{4.14}) and (\ref{4.15}) we also get:
\begin{equation} \label{4.27}
V_n^{\pm}(t)W_{\a, \Lambda}(n,t,\pm i)=2W_{\a, \Lambda}(n,t,\pm
i)P_{\pm}.
\end{equation}

Let us differentiate now $W_{\a, \Lambda}$. For this purpose
notice that
\begin{equation}
 \label{4.28}
(\a \pm i I_N)^{-1}(\a - \lambda I_N)^{-1}= \frac{(\a - \lambda
I_N)^{-1}-(\a \pm i I_N)^{-1}}{ \lambda \pm i}.
\end{equation}
Using (\ref{4.6}), (\ref{4.26}) and (\ref{4.28}) we derive
\begin{eqnarray} \nonumber
\frac{d}{dt}&&W_{\a, \Lambda}(n,t,\lambda)= V_n^+(t) \frac{W_{\a,
\Lambda}(n,t,\lambda)-W_{\a, \Lambda}(n,t,i)}{\lambda - i}+
V_n^-(t) \times
\\ \nonumber && \times
 \frac{W_{\a, \Lambda}(n,t,\lambda)-W_{\a, \Lambda}(n,t,-i)}{\lambda +
i}-2\frac{W_{\a, \Lambda}(n,t,\lambda)-W_{\a,
\Lambda}(n,t,i)}{\lambda - i}P_+-
\\
 \label{4.29} &&-
2\frac{W_{\a, \Lambda}(n,t,\lambda)-W_{\a,
\Lambda}(n,t,-i)}{\lambda + i}P_-.
\end{eqnarray}
In view of (\ref{4.27}) we can rewrite (\ref{4.29}) as
\begin{equation} \label{4.30}
\frac{d}{dt}W_{\a, \Lambda}(n,t,\lambda)=F_n(t, \lambda)W_{\a,
\Lambda}(n,t,\lambda) -W_{\a, \Lambda}(n,t,\lambda) \wh F(t,
\lambda),
\end{equation}
where $F_n$ is given by the second relation in (\ref{0.5.-1}) and
\[
\wh F=2((\lambda - i)^{-1}P_++ (\lambda + i)^{-1}P_-).
\]
Thus in view of Theorem \ref{Tm1.1} and formula (\ref{4.30}) the
non-degenerate matrix functions $\{ \wh W_n \}$ given by
\[ \wh W_n(t, \lambda)= \lambda^{-n}
W_{\a, \Lambda}(n,t,\lambda)
 \left[
\begin{array}{lr} (\lambda - i)^n e^{2t(\lambda - i)^{-1}}
 & 0 \\ 0 & (\lambda + i)^n e^{2t(\lambda + i)^{-1}}
\end{array}
\right]
\]
satisfy equations (\ref{0.5.0}), i.e., the compatibility condition
(\ref{0.5.-2}) is valid. As equation (\ref{0.5.-2}) is equivalent
to (\ref{0.5.1}),  the theorem is proved. $\Box$

Theorem \ref{Tm4.1} together with  Corollary \ref{Cy3.4} yields
the following result.
\begin{Cy} \label{Cy4.2} Under the conditions of Theorem
\ref{Tm4.1} the evolution of the Weyl function $\varphi$ of the
system $W_{n+1}(t, \lambda)=G_n(t, \lambda)W_n(t, \lambda)$ is
given by the formula
\begin{equation} \label{4.31}
\varphi (t, \lambda)=i \t_1^* e^{-2t(\a^* +i
I_N)^{-1}}\S_0(t)^{-1}(\lambda I_N - \wt \b(t))^{-1} e^{-2t(\a +i
I_N)^{-1}} \t_2,
\end{equation}
\[
\wt \b(t) =\a -i e^{-2t(\a +i I_N)^{-1}} \t_2 \t_2^* e^{-2t(\a^*
-i I_N)^{-1}} \S_0(t)^{-1}.
\]
\end{Cy}
As an illustration let us consider a simple example.
\begin{Ee} \label{Ex}
Put $m=n=1$ and $\a =i h$ $(h>0, \, h \not= 1 )$, and choose
scalars $\t_1$, $\t_2$ such that  $|\t_1|^2+ |\t_2|^2=2 h$.
 Then $\a$, $\t_1$, $\t_2$ form an admissible triple and
$\a$, $\S_0(0)=1$, $\Lambda_0(0)=[\t_1 \quad \t_2]$ satisfy the
conditions of Theorem \ref{Tm4.1} and Corollary \ref{Cy4.2}.
Therefore by $($\ref{4.5}$)$ we have
\begin{equation} \label{4.32}
\Lambda_n(t)=h^{-n} \Big[ (h+1)^n \t_1 \exp \Big\{ \frac{2 i
t}{h-1} \Big\} \quad (h-1)^n \t_2 \exp \Big\{ \frac{2 i t}{h+1}
\Big\} \Big].
\end{equation}
From $($\ref{1.1}$)$ and $($\ref{4.32}$)$ it follows that
\begin{equation} \label{4.33} \displaystyle
\S_n(t) \equiv  \frac{c_n(h)}{2h^{2n+1}}, \quad c_n(h):=(h+1)^{2n}
|\t_1|^2 + (h-1)^{2n} |\t_2|^2.
\end{equation}
According to $($\ref{0.6}$)$, $($\ref{4.32}$)$, and
$($\ref{4.33}$)$ we get now
\[
(S_n(t))_{11}= 1- \frac{8h^2|\t_1 \t_2|^2
(h^2-1)^{2n}}{c_n(h)c_{n+1}(h)}, \quad
(S_n(t))_{22}=-(S_n(t))_{11},
\]
\begin{eqnarray} \nonumber
(S_n(t))_{12}= (S_n(t))_{21}^*&=& \frac{4 h \overline{\t_1} \t_2}
{c_n(h)c_{n+1}(h)}  \exp \Big\{ \frac{4 i t}{1-h^2} \Big\} \times
\\ \nonumber
&& \times (h^2-1)^n \big((h+1)^{2n+1} |\t_1|^2 - (h-1)^{2n+1}
|\t_2|^2 \big) .
\end{eqnarray}
Finally Corollary \ref{Cy4.2} yields: \[ \displaystyle \varphi(t,
\lambda)= \exp \Big\{ \frac{4 i t}{1-h^2} \Big\} \frac{i
\overline{\t_1} \t_2 }{ \lambda+i(|\t_2|^2-h)}.
\]
\end{Ee}

\vspace{3em}

M.A. Kaashoek \\ Department of Mathematics, Vrije Universiteit
Amsterdam,  De Boelelaan 1081a, 1081 HV Amsterdam, The
Netherlands. \\

A.L. Sakhnovich \\ Branch of Hydroacoustics,  Marine Institute of
Hydrophysics, National Academy of Sciences, Preobrazhenskaya 3,
Odessa, Ukraine.
\end{document}